\documentclass[11pt]{article}
\usepackage{amsthm}
\usepackage{amsmath}
\usepackage{amsfonts}
\usepackage{cancel}
\usepackage{amssymb}

\usepackage[all, knot]{xy}

\newcommand{\Dorfman}[1]{\bleft  #1\bright}

\newcommand {\emptycomment}[1]{}

\newcommand{\dM}{\mathrm{d}}

\def\C{\mathbb{C}}

\newcommand{\cale}{{\cal E}}
\newcommand{\calf}{{\cal F}}

\newcommand{\br}[1]{[\cdot,\cdot]}

\newcommand{\maps}{\colon}

\newcommand{\g}{R}
\newcommand{\h}{H}

\newcommand{\calC}{{\cal C}}

\newcommand{\calL}{{\cal L}}

\newcommand{\calV}{{\cal V}}

\newcommand{\id}{\mbox{id}}

\newcommand{\ad}{\operatorname{ad}}

\newcommand{\Inn}{\operatorname{Inn}}
\newcommand{\Der}{\operatorname{Der}}

\newcommand{\Cur}{\operatorname{Cur}}

\newcommand{\Ker}{\operatorname{Ker}}

\newcommand{\cendM}{{\operatorname{Cend}(M)}}
\newcommand{\cendr}{{\operatorname{Cend}(R)}}

\newcommand{\half}{\textstyle{\frac{1}{2}}}

\newcommand{\four}{\textstyle{\frac{1}{4}}}

\newcommand{\eight}{\textstyle{\frac{1}{8}}}

\newcommand{\bleft}{[\![}
\newcommand{\bright}{]\!]}

\newcommand{\lam}{\lambda}
\newcommand{\ulam}{_\lambda}
\newcommand{\umu}{_\mu}
\newcommand{\ulamu}{_{\lambda+\mu}}

\newcommand{\ulami}{_{\lambda_1}}
\newcommand{\ulamii}{_{\lambda_2}}
\newcommand{\ulamiii}{_{\lambda_3}}
\newcommand{\ulamipii}{_{\lambda_1+\lambda_2}}
\newcommand{\ulamipiii}{_{\lambda_1+\lambda_3}}
\newcommand{\ulamiipiii}{_{\lambda_2+\lambda_3}}
\newcommand{\ulamipiipiii}{_{\lambda_1+\lambda_2+\lambda_3}}
\newcommand{\ulamIII}{_{\lambda_1+\lambda_2+\lambda_3}}

\newcommand{\uplam}{_{-\partial-\lambda}}

\newcommand{\upmu}{_{-\partial-\mu}}

\newcommand{\om}{\omega}

\def\op{{\oplus}}

\def\la{\langle}
\def\ra{\rangle}

\newcommand{\trl}{\triangleleft}
\newcommand{\trr}{\triangleright}

\newcommand{\pf}{\noindent{\bf Proof.}\ }
\allowdisplaybreaks

\newtheorem{Theorem}{Theorem}[section]
\newtheorem{Proposition}[Theorem]{Proposition}
\newtheorem{Definition}[Theorem]{Definition}

\newtheorem{Example}[Theorem]{Example}
\newtheorem{Remark}[Theorem]{Remark}

\title{Conformal Lie 2-algebras and conformal omni-Lie algebras}

\author{Tao Zhang}

\date{}

\begin{document}



 \maketitle

 \setcounter{section}{0}

 \vskip0.1cm

{\bf Abstract}\quad
The notions of conformal Lie 2-algebras and conformal omni-Lie algebras are introduced.
It is proved that the category of conformal Lie 2-algebras and the category of 2-term conformal $L_{\infty}$-algebras are equivalent.
We construct conformal Lie 2-algebras from conformal omni-Lie algebras and Leibniz conformal algebras.


\section{Introduction}
The algebraic theory of Lie 2-algebras, which is a categorification of Lie algebras, was extensively studied by Baez and Crans \cite{Baez}. As a truncated version of $L_{\infty}$-algebras, Lie 2-algebras replace the underlying vector space with a 2-vector space and the Jacobi identity with a natural transformation known as the Jacobiator, which must satisfy certain coherence laws. Baez and Crans demonstrated the equivalence between the category of Lie 2-algebras and the category of 2-term $L_{\infty}$-algebras. For further advancements in this area of algebraic structures, refer to \cite{Lang,LST,SL,SL0}.

Omni-Lie algebras, introduced by Weinstein in \cite{Wei}, are the most important examples of Lie 2-algebras. They serve as a linearization of the Courant algebroid \cite{LWX} and can be considered a Lie 2-algebra since every Courant algebroid gives rise to a Lie 2-algebra. The study of omni-Lie algebras extends to various aspects, including their generalization to omni-Lie algebroids and omni-Lie 2-algebras as explored in \cite{CL, KW, SLZ}. In a recent paper \cite{ZL}, the Dirac structures of omni-Lie superalgebras are examined.

In this paper, we address the question of whether there exists a categorification of a Lie conformal algebra or a vertex Lie algebra. We give a positive answer to this question by introducing the concept of a conformal Lie 2-algebra, which is a Lie conformal algebra in the category of $\C[\partial]$-modules of 2-vector spaces. Additionally, we define 2-term conformal $L_{\infty}$-algebras and establish the equivalence between the category of conformal Lie 2-algebras and the category of 2-term conformal $L_{\infty}$-algebras.

The second part of this paper is devoted to constructing examples of conformal Lie 2-algebras.
In the third subsection of Section 4, we introduce the concept of conformal omni-Lie algebras and construct conformal Lie 2-algebras from them. Furthermore, in the last part of Section 4, we provide a method to construct conformal Lie 2-algebras from Leibniz conformal algebras.
To clarify this point, we also provide examples since not every Leibniz conformal algebra comes from an omni-Lie algebra.

The organization of this paper is as follows.
In Section 2, we recall some notations and facts about Lie conformal algebras and Leibniz conformal algebras.
In Section 3, we  introduce the notion of conformal Lie 2-algebras and 2-term conformal $L_\infty$-algebras.
It is proved that there is an equivalence between the category of conformal Lie 2-algebras and the category of 2-term conformal $L_{\infty}$-algebras.
In Section 4, we investigate some special cases of conformal Lie 2-algebras, such as skeletal and strict ones.
We also define conformal omni-Lie algebra $\mathcal{E}$ for a Lie conformal algebra and its representation.
At last, we construct conformal Lie 2-algebras from  conformal omni-Lie algebras and Leibniz conformal algebras.

Throughout this paper, we denote by $\mathbb{C}$ the field of complex numbers.
Let $M$ be a vector space.  The space of polynomials of $\lam$ with coefficients in $M$ is denoted by $M[\lam]$.
\section{Preliminaries}

In this section, we will recall some facts and definitions about Lie conformal algebras and conformal Leibinz algebras, see \cite{BKV,DK,GT,HB,Kac,Li,Lib,zhangjiao} for more details.


An \emph{conformal algebra} $R$ is a $\mathbb{C}[\partial]$-module with a $\lambda$-product $\cdot_\lambda \cdot$ which defines a $\mathbb{C}$-linear
map from $A\otimes A\rightarrow A[\lambda]$ satisfying the conformal sesquilinearity condition
\begin{eqnarray}
&&(\partial x)_\lambda y=-\lambda x_\lambda y,~~~x_\lambda \partial y=(\partial+\lambda)x_\lambda y,
\end{eqnarray}
for all $x, y\in A$

\begin{Definition}
A  \emph{Lie conformal algebra} $R$ is a conformal algebra with the $\mathbb{C}$-bilinear
map $[\cdot_\lambda \cdot]: R\times R\rightarrow  R[\lambda]$ satisfying
\begin{eqnarray}
&&[x_\lambda y]=-[y_{-\partial-\lambda}x],~~~~\text{(skew-symmetry)}\\
&&[x_\lambda[y_\mu z]]=[[x_\lambda y]_{\lambda+\mu} z]+[y_\mu[x_\lambda z]],~~~~~~\text{(Jacobi identity)}
\end{eqnarray}
for $x, y, z\in R$.
\end{Definition}

Using the skew-symmetry condition, we can rewrite the Jacobi identity as the following identity:
\begin{eqnarray}\label{eqjacobi}
&&[[x\ulam y]\ulamu z]= [x\ulam[y\umu z]] +  [[x\ulam z]\upmu y].
\end{eqnarray}
\begin{Definition}\label{def:lLeibniz}
A left  \emph{Leibniz conformal algebra} $(R,[\cdot_{\lambda}\cdot])$ is a conformal algebra satisfying the following left Leibniz identity
\begin{eqnarray}
x\circ_\lambda (y\circ_\mu z)=(x\circ_\lambda y)\circ_{\lambda+\mu} z+y\circ_\mu (x\circ_\lambda z).
\end{eqnarray}
\end{Definition}

\begin{Definition}\label{def:rLeibniz}
A right  \emph{Leibniz conformal algebra} $(R,[\cdot_{\lambda}\cdot])$ is a conformal algebra satisfying the following left Leibniz identity
\begin{eqnarray}
(x\circ_\lam y)\circ_{\lambda+\mu}z=x\circ_\lam(y\circ_\mu z)+(x\circ_\lam z)\circ_{\upmu}y.
\end{eqnarray}
\end{Definition}

\begin{Remark}
Using the equivalent identity \eqref{eqjacobi}, we can see a {Lie conformal algebra} as both a left and a right Leibniz conformal algebra.
 However, it should be noted that  the converse is not true as the general Leibniz conformal algebra does not meet the requirement of skew-symmetry.
\end{Remark}

\begin{Definition}
A \emph{left module} $M$ over a Lie conformal algebra $R$ is a $\mathbb{C}[\partial]$-module endowed with a map
$R\times M\longrightarrow M[\lambda]$, $(x, v)\mapsto x\trr_\lambda v$, satisfying the following axioms:
\begin{eqnarray}
&&(\partial x)\trr_\lambda v=-\lambda x\trr_\lambda v,~~~x\trr_\lambda(\partial v)=(\partial+\lambda)x\trr_\lambda v,\\
&& [x_\lambda y]\trr_{\lambda+\mu}v=x\trr_\lambda(y\trr_\mu v)-y\trr_\mu(x\trr_\lambda v).
\end{eqnarray}
for all $x,y\in R, v\in M$.
Similarly,  a \emph{right module} $M$ over a Lie conformal algebra $R$ is a $\mathbb{C}[\partial]$-module endowed with a bilinear map
$M\times R\longrightarrow M[\lambda]$, $(v, x)\mapsto v\trl_\lambda x$, satisfying the following axioms:
\begin{eqnarray}
&&(\partial v)\trl_\lambda x=-\lambda v\trl_\lambda x,~~v\trl_\lambda (\partial x)=(\partial+\lambda)v\trl_\lambda x,\\
&&v\trl_\mu[x_\lambda y]=(v\trl_\mu x)\trl_{\lambda+\mu}y-(v\trl_\mu y)\trl_{\uplam}x.
\end{eqnarray}
A right $R$-module become  a left $R$-module if we define $x\trr_\lambda v$ to be $-v\trl\uplam x$.
\end{Definition}

\begin{Definition}
A  \emph{module} $M$ over a Leibniz conformal algebra $(R, [\cdot_\lambda\cdot])$ is a $\mathbb{C}[\partial]$-module endowed with two $\mathbb{C} $-bilinear maps
$R\times M\longrightarrow M[\lambda]$, $(x, v)\mapsto x\trr_\lambda v$ and $M\times R\longrightarrow M[\lambda]$, $(v, x)\mapsto v{\trl}_\lambda x$ satisfying the following axioms $(a, b\in R, v\in M)$:
\begin{eqnarray}
 (\partial x)\trr_\lambda v&=&-\lambda x\trr_\lambda v,~~~x\trr_\lambda(\partial v)=(\partial+\lambda)x\trr_\lambda v,\\
 (\partial v)\trl_\lambda x&=&-\lambda v\trl_\lambda x,~~v\trl_\lambda (\partial x)=(\partial+\lambda)v\trl_\lambda x,\\
x{\trr}_\lambda (y{\trr}_\mu v)&=&(x\circ_\lambda y){\trr}_{\lambda+\mu}v
+y{\trr}_\mu(x{\trr}_\lambda v),\\
x{\trr}_\lambda(v{\trl}_\mu y)&=&(x{\trr}_\lambda v){\trl}_{\lambda+\mu}y+v{\trl}_\mu[x_\lambda y],\\
v{\trl}_\lambda(x\circ_\mu y)&=&(x{\trl}_\lambda x){\trl}_{\lambda+\mu}y+x{\trr}_\mu(v{\trl}_\lambda y).
\end{eqnarray}
We denote it by $(M,{\trr}_\lambda, {\trl}_\lambda)$.
\end{Definition}

Next, we give some examples of Lie conformal algebras.
\begin{Example}
Let $\mathfrak{g}$ be a Lie algebra. Let $\Cur\mathfrak{g}:=\mathbb{C}[\partial] \otimes \mathfrak{g}$ be the free $\mathbb{C}[\partial]$-module. Then $\Cur \mathfrak{g}$ is a Lie conformal algebra, called current Lie conformal algebra, with $\lambda$-bracket given by:
$$[(f(\partial) \otimes x)_\lambda(g(\partial) \otimes y)]:=f(-\lambda) g(\lambda+\partial) \otimes[x, y],$$
for all $ f(\partial), g(\partial) \in \mathbb{C}[\partial], x, y\in \mathfrak{g}$.

The Virasoro Lie conformal algebra $\text{Vir}$ is the simplest nontrivial
example of Lie conformal algebras. It is defined by
$$\text{Vir}=\mathbb{C}[\partial]L, ~~[L_\lambda L]=(\partial+2\lambda)L.$$
Coeff$\text{(Vir)}$ is just the Witt algebra.
\end{Example}

 Let $M$ and $N$ be $\mathbb{C}[\partial]$-modules. A conformal linear map $f$ from $M$ to $N$ is a $\mathbb{C}$-linear map
 $f_\lambda: M \rightarrow \mathbb{C}[\lambda] \otimes N$ such that $f_\lambda \partial=(\partial+\lambda) f_\lambda$.
The category of  $\mathbb{C}[\partial]$-modules with conformal linear maps  as morphisms is denoted by ${\rm\mathbf{Vect}}^{\partial}$.

Let $\operatorname{Chom}(M, N)$ denote the set of conformal linear maps from $M$ to $N$. Then $\operatorname{Chom}(M, N)$ is a $\mathbb{C}[\partial]$-module via:
$$
\partial f_\lambda=-\lambda f_\lambda.
$$
The composition $f g: L \rightarrow N$ of conformal linear maps $f: M \rightarrow N$ and $g: L \rightarrow M$ is given by
$\left(f_\lambda g\right)_{\lambda+\mu}=f_\lambda g_\mu$.
If $M$ is a finitely generated $\mathbb{C}[\partial]$-module, then $\operatorname{Cend}(M):=\operatorname{Chom}(M, M)$ is an associative conformal algebra with respect to the above composition. Thus, $\operatorname{Cend}(M)$ becomes a Lie conformal algebra, which is denoted as $\operatorname{Cend}(M)$, with respect to the following $\lambda$-bracket
$$\left[f_\lambda g\right]_\mu=f_\lambda g_{\mu-\lambda}-g_{\mu-\lambda} f_\lambda,$$
equivalently,
$$\left[f_\lambda g\right]=f_\lambda g-g_{-\partial-\lambda} f.$$
All $\mathbb{C}[\partial]$-modules  in this paper are assumed to be finitely generated.

A homomorphism  between two Lie conformal algebras $({R},[\cdot\ulam\cdot ])$ and $({R}',[\cdot\ulam\cdot]')$ is a conformal
linear map ${f}: {R} \to {R}'$ such that
$$f_\lambda \partial(x)=(\partial+\lambda) f_\lambda(x),\quad {f}([x_\lambda y]) = [{f}(x)_\lambda {f}(y)]' $$ 
for all $x, y\in  {R}$.

An $n$-cochain ($n\in \mathbb{Z}_{\geq0}$) of a Lie conformal algebra $R$ with coefficients
 in a module $M$ is a $\mathbb{C}$-linear map
\begin{eqnarray*}
\gamma:R^{\otimes n}\rightarrow M[\lambda_1,... ,\lambda_{n-1}],~~(x_1,...,x_n)\mapsto \gamma_{\lambda_1,...,\lambda_{n-1}}(x_1,... , x_n)=\gamma(x_1{}_{\lam_1}\ldots x_{n-1}{}_{\lam_{n-1}} x_n),
\end{eqnarray*}
where $M[\lambda_1,...,\lambda_{n-1}]$ denotes the space of polynomials with coefficients in $M$, satisfying the following conditions:

(1) Conformal linearity:
\begin{eqnarray*}
\gamma_{\lambda_1,...,\lambda_{n-1}}(x_1,... ,\partial x_i,... ,x_n)=-\lambda_i\gamma_{\lambda_1,...,\lambda_{n-1}}(x_1,... ,x_i,... ,x_n).
\end{eqnarray*}

(2) Skew-symmetry:
\begin{eqnarray*}
&&\gamma(x_{1},..., x_{i+1}, \partial(x_{i}),..., x_{n})=-\gamma(x_{1}, ..., x_{i}, \partial(x_{i+1}), ... , x_{n}).
\end{eqnarray*}

Let $R^{\otimes 0} = \mathbb{C}$ as usual so that a $0$-cochain  is an element of $a\in M$ and $(\delta\gamma)_\lambda x=x_\lambda \gamma$. Define a differential $\delta$ of a cochain $\g$ by
\begin{eqnarray*}
&&(\delta\gamma)_{\lambda_1,...,\lambda_{n}}(x_{1},...,x_{n+1})\nonumber\\
&=&\sum_{i=1}^{ n+1}(-1)^{i+1}{x_{i}}\trr_{\lambda_i}\gamma_{\lambda_1,...,\hat{\lambda}_{i},...,\lambda_{n}}(x_{1},..., \hat{x}_{i},..., x_{n+1})\\
  &&+\sum_{1\leq i<j\leq n+1} (-1)^{i+j}\gamma_{\lambda_i+\lambda_j, \lambda_1,...,\hat{\lambda}_i,...\hat{\lambda}_j,...,\lambda_{n}}([ x_{i}{}_{\lambda_i}x_{j}], x_{1},... ,\hat{x}_{i},... , \hat{x}_{j},...,x_{n+1}).
\end{eqnarray*}
It is proved that $\delta^2=0$ \cite{BKV}.
Therefore $\{C^*(\g,M), \delta\}$ is indeed a cochain complex, whose cohomology is called the
cohomology of the Lie conformal algebra $\g$ with coefficients in the representation $M$.

Using the skew-symmetry, we can also write the differential $\delta$ as follows:
\begin{eqnarray*}
&&(\delta\gamma)_{\lambda_1,...,\lambda_{n}}(x_{1},...,x_{n+1})\nonumber\\
&=&\sum_{i=1}^{ n}(-1)^{i+1}{x_{i}}\trr_{\lambda_i}\gamma_{\lambda_1,...,\hat{\lambda}_{i},...,\lambda_{n}}(x_{1},..., \hat{x}_{i},..., x_{n+1})\\
&&+ (-1)^{n+1}\gamma_{\lambda_1,...,\lambda_{n-1}} (x_1, \ldots, x_{n}){\trl}_{\lam_1+\ldots+\lam_n} x_{n+1}\\
  &&+\sum_{1\leq i<j\leq n+1} (-1)^{i}\gamma_{ \lambda_1,...,\hat{\lambda}_i,...,\lambda_i+\lambda_j,...,\lambda_{n}}(x_{1},... ,{x}_{j-1},[ x_{i}{}_{\lambda_i}x_{j}], {x}_{j+1},...,x_{n+1}).
\end{eqnarray*}
This formula was first appeared in \cite{zhangjiao} as the cohomology of a Leibniz conformal algebras.

\begin{Remark}
We use a slightly different notation as in \cite{BKV,DK} where they use $\gamma_{\lambda_1,...,\hat{\lambda}_{i},...,\lambda_{n+1}}$ to denote the map  $\gamma$, but in fact the index $\lambda_{n+1}$ always does't appear in the calculation of the cohomology theory of a Lie conformal algebra, so we use $\gamma_{\lambda_1,...,\hat{\lambda}_{i},...,\lambda_{n}}$ instead of it.
\end{Remark}

\section{Conformal Lie 2-algebras and 2-term conformal $L_{\infty}$-algebras}

In this section, we introduced the concept of conformal Lie 2-algebras and 2-term conformal $L_{\infty}$-algebras.
It is proved that the category of conformal Lie $2$-algebras and the category of $2$-term conformal $L_{\infty}$-algebras are equivalent.

\subsection{Conformal Lie 2-algebras}

Let's denote the category of  $\mathbb{C}[\partial]$-modules or $\partial$-vector spaces by ${\rm\mathbf{Vect}}^{\partial}$.

\begin{Definition}
A {\bf conformal 2-vector space} is a category in ${\rm\mathbf{Vect}}^{\partial}$.
\end{Definition}

More precisely, a conformal $2$-vector space $\calV$ is a category with a $\mathbb{C}[\partial]$-module $\calV_0$ and a $\mathbb{C}[\partial]$-module $\calV_1$. The category has source and target maps $s,t \maps \calV_{1} \rightarrow \calV_{0}$, an identity-assigning map $i \maps \calV_{0} \rightarrow \calV_{1}$, and a composition map $\circ \maps \calV_{1} \times_{\calV_{0}} \calV_{1} \rightarrow \calV_{1}$. These maps are all conformal linear maps. A morphism $f$ from source $x$ to target $y$ is denoted by $f \maps x \rightarrow y$, where $s(f) = x$ and $t(f) = y$. The notation $i(x)$ is also written as $1_x$.

Conformal 2-vector spaces are in one-to-one correspondence with 2-term complexes of $\mathbb{C}[\partial]$-modules. A 2-term complex of $\mathbb{C}[\partial]$-modules is a pair of $\mathbb{C}[\partial]$-modules $\calC_1$ and $\calC_0$ with a differential between them: $\calC_1\stackrel{\dM}{\longrightarrow}\calC_0$. Given a conformal 2-vector space $\calV$, $\Ker(s)\stackrel{t}{\longrightarrow}\calV_0$ is a 2-term complex of $\mathbb{C}[\partial]$-modules. Conversely, any 2-term complex of $\mathbb{C}[\partial]$-modules $\calV_1\stackrel{\dM}{\longrightarrow}\calV_0$ gives rise to a conformal 2-vector space. In this conformal 2-vector space, the set of objects is $\calC_0$ and the set of morphisms is $\calC_0\oplus \calC_1$. The source map $s$ is given by $s(x,h)=x$, and the target map $t$ is given by $t(x,h)=x+\dM h$, where $x\in \calV_0$ and $h\in \calV_1$. The conformal 2-vector space associated to the 2-term complex $\calV_1\stackrel{\dM}{\longrightarrow}\calV_0$ is denoted by $\calV$:
\begin{equation}\label{eqn:V}
\calV=\begin{array}{c}
\calV_1:=\calC_0\oplus \calC_1\\
\vcenter{\rlap{s }}~\Big\downarrow\Big\downarrow\vcenter{\rlap{t }}\\
\calV_0:=\calC_0.
 \end{array}\end{equation}

\begin{Definition} \label{defnlie2alg}
A {\bf conformal Lie 2-algebra} consists of a  conformal 2-vector space $\calL$ equipped with
\begin{itemize}
\item a skew-symmetric conformal sesquilinear functor, the {\bf bracket}, $[\cdot\ulam \cdot]\maps \calL \times \calL\rightarrow \calL[\lambda]$
\item a conformal sesquilinear natural isomorphism, the {\bf Jacobiator},
$$J_{x\ulami y\ulamii z} \maps [x\ulami [y\ulamii z]] \to [[x\ulami y]\ulamipii z] + [y\ulamii [x\ulami z]],$$
\end{itemize}
such that the following  {\bf Jacobiator identity} is satisfied
\begin{align}\label{Jacobiator}
&\quad J_{x\ulami y\ulamii [z\ulamiii t]} \Big(1+ [y\ulamii J_{x\ulami z\ulamiii t}]\Big) \Big({J_{[x\ulami y]\ulamipii z\ulamiii t}+ J_{y\ulamii [x\ulami z]\ulamipiii t}+ J_{y\ulamii z\ulamiii [x\ulami t]}}\Big)\notag\\
&=[x\ulami J_{y\ulamii z\ulamiii t}] \Big(J_{x\ulami [y\ulamii z]\ulamiipiii t}+ J_{x\ulami z\ulamiii [y\ulamii t]}\Big) \Big([J_{x\ulami y\ulamii z}{}\ulamipiipiii t]+ 1+1+ [z\ulamiii J_{x\ulami y\ulamii t}] \Big).
\end{align}
\end{Definition}

We represent the Jacobiator identity in the following commutative diagram, which illustrates the relationship between two methods of utilizing the Jacobiator to rebracket the expression $[x\ulami[y\ulamii[z\ulamiii t]]]$:
$$\def\objectstyle{\scriptstyle}
  \def\labelstyle{\scriptstyle}
\xymatrix{
&[x\ulami[y\ulamii[z\ulamiii t]]]\ar[dr]^{ J_{x\ulami y\ulamii [z\ulamiii t]}}\ar[dl]_{[x\ulami J_{y\ulamii z\ulamiii t}] }&\\
[x\ulami [[y\ulamii z]\ulamiipiii t]]+ [x\ulami [z\ulamiii [y\ulamii t]]]\ar[dd]^{J_{x\ulami [y\ulamii z]\ulamiipiii t}+ J_{x\ulami z\ulamiii [y\ulamii t]}}
&& [[x\ulami y]\ulamipii [z\ulamiii t]]+ [y\ulamii [x\ulami [z\ulamiii t]]]\ar[dd]_{1+ [y\ulamii J_{x\ulami z\ulamiii t}]}\\
&&\\
{\begin{aligned}&\scriptstyle [[x\ulami [y\ulamii z]]\ulamipiipiii t]+  [[y\ulamii z]\ulamiipiii[x\ulami t]]\\[-.5em]
&\scriptstyle+ [[x\ulami z]\ulamipiii [y\ulamii t]]+  [z\ulamiii [x\ulami [y\ulamii t]]]\end{aligned}}
\ar[dr]_{[J_{x\ulami y\ulamii z}{}\ulamipiipiii t]+ 1+1+ [z\ulamiii J_{x\ulami y\ulamii t}]\qquad\qquad\quad}&&
{\begin{aligned}&\scriptstyle[[x\ulami y]\ulamipii [z\ulamiii t]]+ [y\ulamii [[x\ulami z]\ulamipiii t]] \\[-.5em]
&\scriptstyle \quad+   [y\ulamii [z\ulamiii [x\ulami t]]] \end{aligned}}
\ar[dl]^{\qquad {J_{[x\ulami y]\ulamipii z\ulamiii t}+ J_{y\ulamii [x\ulami z]\ulamipiii t}+ J_{y\ulamii z\ulamiii [x\ulami t]}}}\\
&P=Q&}
\\ \\
$$
where $P$ and $Q$ are given by
\begin{eqnarray*}
  P &=&[[[x\ulami y]\ulamipii z]\ulamipiipiii t]+ [[y\ulamii [x\ulami z]]\ulamipiipiii t]\\
  &&+ [[y\ulamii z]\ulamiipiii[x\ulami t]]+ [[x\ulami z]\ulamipiii [y\ulamii t]]\\
   &&+  [z\ulamiii [[x\ulami y]\ulamipii t]]+   [z\ulamiii [y\ulamii [x\ulami t]]]\qquad=Q.
\end{eqnarray*}

If we drop the skew-symmetry, we obtain the concept of  conformal Leibinz 2-algebras.
\begin{Definition} \label{Leib2alg}
A {\bf conformal Leibinz 2-algebra} consists of a  conformal 2-vector space $\calL$ equipped with
 a conformal sesquilinear functor $[\cdot\ulam \cdot]\maps \calL \times \calL\rightarrow \calL[\lambda]$
and a conformal sesquilinear natural isomorphism
$$J_{x\ulami y\ulamii z} \maps [x\ulami [y\ulamii z]] \to [[x\ulami y]\ulamipii z] + [y\ulamii [x\ulami z]],$$
such that the above \eqref{Jacobiator} is satisfied.
\end{Definition}

In fact, all the results in the following of this section can be easily generalized to the realm of conformal Leibinz 2-algebra, but we omit the details.

\begin{Definition}
  Given two conformal Lie 2-algebras $(\calL,[\cdot\ulam\cdot], J)$ and $(\calL^\prime,[\cdot\ulam\cdot]', J^\prime)$, a conformal Lie 2-algebra morphism $F:\calL\longrightarrow
  \calL^\prime$ consists of:
  \begin{itemize}
   \item[$\bullet$] a conformal  functor $(F^0,F^1)$ from the underlying conformal  2-vector space of $\calL$ to that of $\calL^\prime$;

 \item[$\bullet$] a  natural transformation
 $$
F^2\ulam(x,y):[F^0(x)\ulam F^0(y)]'\longrightarrow F^0([x\ulam y])
 $$
such that the following diagram commutes:
  \end{itemize}
 $$
\footnotesize{ \xymatrix{
[F^0(x)\ulam [F^0(y)\umu F^0(z)]']'\ar[d]_{J_{F^0(x)\ulam F^0(y)\umu F^0(z)}}\ar[rr]^{\qquad\qquad\qquad[1,F^2\umu(y, z)]'
\quad\qquad\quad}&&[F^0(x)\ulam F^0([y\umu z])]'\ar[d]^{F^2\ulam (x, [y\umu z])}\\
~[[F^0(x)\ulam F^0(y)]'\ulamu F^0(z)]'+ [F^0(y)\umu [F^0(x)\ulam F^0(z)]']'\ar[d]_{[F^2\ulam (x,y),1]'+[1,F^2\ulam (x,z)]'}&&F^0([x\ulam [y\umu z]])\ar[d]^{F^0 J_{x\ulam y\umu z}}\\
~[F^0([x\ulam y])\ulamu F^0(z)]'+ [F^0(y)\umu F^0[x\ulam z]]'\ar[rr]^{\footnotesize\qquad\quad {F^2\ulamu ([x\ulam y],z)+F^2\umu (y, [x\ulam z])}}&&
\footnotesize{F^0([[x\ulam y]\ulamu z])+ F^0([y\umu [x\ulam z]]). }}}
$$
\end{Definition}

The identity morphism ${\id}_\calL:\calL\longrightarrow \calL$ is associated with the identity functor and an identity natural transformation $({\id}_\calL)_2$.
The composition of a pair of conformal Lie 2-algebra morphisms, denoted by $G\circ F$, is defined as follows: let $\calL,~\calL'$ and $\calL''$ be conformal Lie 2-algebras. The functor $((G\circ F)^0,(G\circ F)^1)$ is obtained by composing $(G^0,G^1)$ and $(F^0,F^1)$, while $(G\circ F)^2$ is defined as the following composition:
$$
 \xymatrix{
[G^0\circ F^0(x)\ulam G^0\circ F^0(y)]''\ar[dd]_{G^2\ulam(F^0(x),F^0(y))}\ar[dr]^{(G\circ F)^2\ulam(x,y)}&&\\
&G^0\circ F^0[x\ulam y].&\\
 G^0[F^0(x)\ulam F^0(y)]'\ar[ur]_{G^1(F^2\ulam(x,y))}&&,}.
$$
It is easy to see that there is a category {\bf LieC2Alg} with conformal Lie 2-algebras as objects and conformal Lie 2-algebra morphisms as morphisms.

\subsection{2-term conformal $L_{\infty}$-algebras}

Now we introduced the concept of  2-term conformal $L_{\infty}$-algebras.
For general theory of $L_{\infty}$-algebras with applications, see \cite{B2,LS93,Kon,Roy,SS}
\begin{Definition} \label{2termliealgebra}
A $2$-term conformal $L_{\infty}$-algebra $\calV=\calV_0\oplus \calV_1$ is a complex consisting of the following
data:
\begin{itemize}
  \item two $\mathbb{C}[\partial]$-modules $\calV_{0}$ and   $\calV_{1}$ together with a conformal linear map  $\dM\maps V_{1} \rightarrow V_{0}$.

  \item a conformal sesquilinear map $l^2_{\lam}\maps \calV_{i} \times \calV_{j}   \rightarrow \calV_{i+j}[\lam],$ where $0 \leq i + j \leq 1$,

  \item a conformal sesquilinear map $l^3_{\lam_1,\lam_2}\maps \calV_{0} \times \calV_{0} \times \calV_{0} \rightarrow \calV_{1}[\lam_1,\lam_2]$.
\end{itemize}
These maps satisfy the following conditions:
\begin{itemize}
  \item[(a)] $l^2_{\lam}(x,y) + l^2_{\uplam}(y,x)=0$,
  \item[(b)] $l^2_{\lam}(x,h) +l^2_{\uplam}(h,x)=0$,
  \item[(c)] $l^2_{\lam}(h,k)=0$,
  \item[(d)] $l^3_{\lam_1,\lam_2}(x,y,z)$ is totally skew symmetric,
  \item[(e)] $\dM(l^2_{\lam}(x,h)) = l^2_{\lam}(x,\dM h)$,
  \item[(f)] $l^2_{\lam}(\dM h,k) =  l^2_{\lam}(h,\dM k)$,
  \item[(g)] $\dM(l^3_{\lam_1,\lam_2}(x,y,z))=l^2_{\lam_1}(x,l^2_{\lam_2}(y,z))-l^2_{\lam_1+\lam_2}(l^2_{\lam_1}(x,y),z)- l^2_{\lam_2}(y,l^2_{\lam_1}(x,z))$,
  \item[(h)] $l^3_{\lam_1,\lam_2}(x,y,\dM h)=l^2_{\lam_1}(x,l^2_{\lam_2}(y,h))-l^2_{\lam_1+\lam_2}(l^2_{\lam_1}(x,y),h)-l^2_{\lam_2}(y,l^2_{\lam_1}(x,h))$,
  \item[(i)]
    \begin{flalign*}
    &\delta l^3_{\lam_1,\lam_2,\lam_3}(x,y,z,t)\\
    =&l^2_{\lam_1}(x, l^3_{\lam_2,\lam_3}(y, z,t)) -   l^2_{\lam_2}(y, l^3_{\lam_1,\lam_3}(x, z,t))+ l^2_{\lam_3}(z, l^3_{\lam_1,\lam_2}(x, y,t))\\
     &+l^2_{\lam_1+\lam_2+\lam_3}(l^3_{\lam_1,\lam_2}(x, y, z),t)- l^3_{\lam_1+\lam_2,\lam_3}(l^2_{\lam_1}(x, y), z,t)\\
     &-  l^3_{\lam_2,\lam_1+\lam_3}(y,l^2_{\lam_1}(x, z),t)- l^3_{\lam_2,\lam_3}(y, z,l^2_{\lam_1}(x,t))+ l^3_{\lam_1,\lam_2+\lam_3}(x,l^2_{\lam_2}(y, z), t)\\
     & + l^3_{\lam_1,\lam_3}(x,z,l^2_{\lam_2}(y,t))- l^3_{\lam_1,\lam_2}(x, y,l^2_{\lam_3}(z,t))\\
     =&0,
      \end{flalign*}
\end{itemize}
for all $x,y,z,t\in \calV_{0}$ and $h, k \in \calV_{1}.$
\end{Definition}


\begin{Definition}\label{defi:Lie-2hom}
Let $(\calV;\dM,l^2_{\lam},l^3_{\lam_1,\lam_2})$ and $(\calV';\dM',{l}^2_{\lam},{l'}^3_{\lam_1,\lam_2})$ be two $2$-term conformal $L_{\infty}$-algebras.
A $CL_\infty$-morphism $f=(f^0,f^1,f^2)$ from $\calV$ to $\calV'$ consists of conformal
 linear maps $f^0:\calV_0\rightarrow \calV_0',~f^1:\calV_{1}\rightarrow \calV_{1}'$
 and $f^2: \calV_{0}\times \calV_0\rightarrow \calV_{1}'[\lam]$,
such that the following equalities hold for all $ x,y,z\in \calV_{0},
a\in \calV_{1},$
\begin{itemize}
\item [$\rm(i)$] $f^0\dM=\dM'f^1$,
\item[$\rm(ii)$] $f^0l^2_{\lam}(x,y)-{l'}^2_{\lam}(f^0(x),f^0(y))=\dM'f\ulam^2(x,y),$
\item[$\rm(iii)$] $f^1l^2_{\lam}(x,a)-{l'}^2_{\lam}(f^0(x),f^1(a))=f\ulam^2(x,\dM a)$,
\item[$\rm(iv)$] $f^1(l^3_{\lam_1,\lam_2}(x,y,z))-{l'}^3_{\lam_1,\lam_2}(f^0(x),f^0(y),f^0(z))$\\
  $=f\ulami^2(x, l^2\ulamii(y,z))- f^2_{\lam_1+\lam_2}(l^2_{\lam_1}(x,y),z) - f^2\ulamii(y,l^2_{\lam_1}(x,z))$\\
   $+ {l'}^2_{\lam_1}(f^0(x), f^2\ulamii(y,z)) - {l'}^2_{\lam_1+\lam_2}(f^2\ulami(x,y), f^0(z))- {l'}^2_{\lam_2}(f^0(y), f^2\ulami(x,z)).$
\end{itemize}
 If $f^2=0$, the $CL_\infty$-morphisms $f$ is called a strict $CL_\infty$-morphisms.
\end{Definition}

Let $f:\calV\to \calV'$ and $g:\calV'\to \calV''$ be two $CL_\infty$-morphisms, then their composition $g\circ f:\calV\to \calV''$
is a $CL_\infty$-morphism defined as $(g\circ f)^0=g^0\circ f^0$, $(g\circ f)^1=g^1\circ f^1$
and
$$(g\circ f)^2\ulam(x,y)=g^2\ulam (f^0(x)f^0(y))+g^1(f^2\ulam(x,y)).$$

The identity $CL_\infty$-morphism $1_\calV: \calV\to \calV$ has the identity chain map together with $(1_\calV)_2=0$.

There is a category {\bf 2CL$_\infty$} with 2-term conformal $L_\infty$-algebras as objects and $CL_\infty$-morphisms as morphisms.

\subsection{Equivalence}

Now we establish the equivalence between the category of conformal Lie $2$-algebras and $2$-term conformal $L_{\infty}$-algebras.

\begin{Theorem} The categories  ${\bf2CL_\infty}$ and  ${\bf LieC2Alg}$ are equivalent.
\end{Theorem}

\pf  First, we show how to construct a conformal Lie 2-algebra from a 2-term conformal $L_\infty$-algebra.

Let $\calV=(\calV_1\stackrel{\dM}{\longrightarrow}\calV_0,l^2_{\lam},l^3_{\lam_1,\lam_2})$ be a 2-term conformal $L_\infty$-algebra, we consider the  conformal  $2$-vector space $\calL=(\calV_0\oplus \calV_1\rightrightarrows \calV_0)$ given by \eqref{eqn:V},
that is $\calL$ has objects $\calL_0=\calV_0$, $\calL_1=\calV_0\oplus \calV_1$ and  morphisms
$f\maps x \rightarrow y$ in $\calL_{1}$ by
$f=(x, h)$ where  $x \in \calV_{0}$ and $h \in \calV_{1}$.  The source, target, and
identity-assigning maps in $\calL$ are given by
\begin{eqnarray*}
  s(f) &=& s(x, h) = x ,\\
  t(f) &=& t(x, h) = x + \dM h, \\
  i(x) &=& (x, 0),
\end{eqnarray*}
and we have $t(f) - s(f) = \dM h$.

Now we define  a  skew-symmetric conformal sesquilinear functor on $\calL$ by
$$[(x,h)\ulam (y,k)]=(l^2_{\lam}(x,y), l^2_{\lam}(x,k)+l^2_{\lam}(h,y)+l^2_{\lam}(\dM h,k))$$
for all $(x,h), (y,k)\in \calV_0\oplus \calV_1$.
The Jacobiator is given as follows
$$J_{x\ulami y\ulamii z}:=([x\ulami [y\ulamii z]], \, l^3_{\lam_1,\lam_2}(x,y,z))$$
where $x,y,z\in \calV_0$.
Then by Condition $(g)$, we have $J_{x\ulami y\ulamii z}$ is a morphism from source $[x\ulami[y\ulamii z]]$ to target $[[x\ulami y]\ulamipii z] +
[y\ulamii [x\ulami z]]$.

Next we show that $J_{x\ulami y\ulamii z}$ is a natural isomorphism.
 We only check naturality in the third variable, the other two cases are similar.  Let
$f \maps z \rightarrow z'.$ Then, $J_{x\ulami y\ulamii z}$ is natural in $z$ if
the following diagram commutes:
$$\xymatrix{
      [x\ulami [y\ulamii z]] \ar[rrr]^{[[1_x,1_y],f]}  \ar[dd]_{J_{x\ulami y\ulamii z}}      &&&   [x\ulami [y\ulamii z']] \ar[dd]^{J_{x\ulami y\ulamii z'}} \\ \\
    [[x\ulami y]\ulamipii z] +   [y\ulamii [x\ulami z]]   \ar[rrr]^{[1_x,[1_y,f]]+ [[1_x,f],1_y]}      &&&  [[x\ulami y]\ulamipii z'] +   [y\ulamii [x\ulami z']]  }$$
\\
\noindent Thus  we only need to show
$$([x\ulami [y\ulamii z]] , l^3_{\lam_1,\lam_2}(x,y,z') + [x\ulami [y\ulamii h]] ) =([x\ulami [y\ulamii z]],   [x\ulami y]\ulamipii h]]+ [y\ulamii [x\ulami h]] +l^3_{\lam_1,\lam_2}(x,y,z)).$$
This holds by condition $(h)$ in Definition \ref{2termliealgebra}.

From condition $(i)$ in Definition \ref{2termliealgebra}, we have
    \begin{flalign*}
    &l^2_{\lam_1}(x, l^3_{\lam_2,\lam_3}(y, z,t)) -   l^2_{\lam_2}(y, l^3_{\lam_1,\lam_3}(x, z,t))+ l^2_{\lam_3}(z, l^3_{\lam_1,\lam_2}(x, y,t))\\
     &+l^2_{\lam_1+\lam_2+\lam_3}(l^3_{\lam_1,\lam_2}(x, y, z),t)- l^3_{\lam_1+\lam_2,\lam_3}(l^2_{\lam_1}(x, y), z,t)\\
     &-  l^3_{\lam_2,\lam_1+\lam_3}(y,l^2_{\lam_1}(x, z),t)- l^3_{\lam_2,\lam_3}(y, z,l^2_{\lam_1}(x,t))+ l^3_{\lam_1,\lam_2+\lam_3}(x,l^2_{\lam_2}(y, z), t)\\
     & + l^3_{\lam_1,\lam_3}(x,z,l^2_{\lam_2}(y,t))- l^3_{\lam_1,\lam_2}(x, y,l^2_{\lam_3}(z,t))=0.
      \end{flalign*}
By the  skew-symmetry of $l^2_{\lam}=[\cdot\ulam\cdot]$, we obtain
\begin{eqnarray*}
     && l^3_{\lam_1,\lam_2}(x, y,[z\ulamiii t])+[y\ulamii l^3_{\lam_1,\lam_3}(x, z,t))]\\
     &&+l^3_{\lam_1+\lam_2,\lam_3}([x\ulami y], z,t) + l^3_{\lam_2,\lam_1+\lam_3}(y,[x\ulami z],t)\\
     &&+ l^3_{\lam_2,\lam_3}(y, z,[x\ulami t])\\
   &=&[x\ulami  l^3_{\lam_2,\lam_3}(y, z,t)] +[z\ulamiii  l^3_{\lam_1,\lam_2}(x, y,t)]\\
     &&+[l^3_{\lam_1,\lam_2}(x, y, z)\ulamIII t)   + l^3_{\lam_1,\lam_2+\lam_3}(x, [y\ulamii z], t)\\
      &&+ l^3_{\lam_1,\lam_3}(x,z,[y\ulamii t]).
\end{eqnarray*}
This is equivalent to Jacobiator identity \eqref{Jacobiator} in Definition \ref{defnlie2alg}.
Thus from a 2-term conformal $L_\infty$-algebra, we obtain a conformal Lie 2-algebra.

For any $CL_\infty$-morphism $f=(f^0,f^1,f^2)$ form $\calV$ to
$\calV'$, we construct a conformal Lie 2-algebra morphism $F=T(f)$
from $L=T(\calV)$ to $L'=T(\calV')$ as follows.

Let $F^0=f^0,~F^1=f^0\oplus f^1$, and $F^2$ be given by
$$
F^2\ulam(x,y)=([f^0(x)\ulam f^0(y)],f^2\ulam(x,y)).
$$
Then $F^2(x,y)$ is a natural isomorphism
from $[F^0(x)\ulam F^0(y)]$ to $F^0[x\ulam y]$, and $F=(F^0,F^1,F^2)$ is a conformal Lie 2-algebra
morphism from $\calL$ to $\calL'$.

One can also deduce that $T$ preserves the identity $CL_\infty$-morphisms and
the composition of $CL_\infty$-morphisms. Thus, $T$ constructed above is a
functor from {\bf 2CL$_\infty$} to {\bf LieC2Alg}.

Conversely, given a conformal Lie 2-algebra $\calL$, we define $l^2_{\lam}$ and
$l^3_{\lam_1,\lam_2}$ on the 2-term complex $\calL_1\supseteq \ker(s)=\calV_1\stackrel{\dM}{\longrightarrow}\calV_0=\calL_0$
by
\begin{itemize}
\item $l_1h = t(h)$ for $h \in \calV_1 \subseteq \calL_1$.
\item $l^2_{\lam}(x,y) = [x\ulam y]$ for $x,y \in \calV_0 = \calL_0$.
\item $l^2_{\lam}(x,h) = [1_x{}\ulam h]$
for $x \in \calV_0 = \calL_0$ and $h \in \calV_1 \subseteq \calL_1$.
\item $l^2_{\lam}(h,k) = 0$ for $h,k \in \calV_1 \subseteq \calL_1$.
\item $l^3_{\lam_1,\lam_2}(x,y,z) =p_1J_{x\ulami y\ulamii z}$ for $x,y,z \in \calV_0 = \calL_0$, where $p_1: \calL_1=\calV_0\oplus \calV_1\longrightarrow \calV_1$ is the projection.
\end{itemize}
Then one can verify that
$(\calV_1\stackrel{\dM}{\longrightarrow}\calV_0,l^2_{\lam},l^3_{\lam_1,\lam_2})$ is a 2-term conformal $L_\infty$-algebra.

Let $F=(F^0,F^1,F^2):\calL\longrightarrow \calL'$ be a conformal Lie 2-algebra
morphism, and $S(\calL)=\calV,~S(\calL')=\calV'$. Define
$S(F)=f=(f^0,f^1,f^2)$ as follows. Let $f^0=F^0$,
$f^1=F^1|_{V_1=\Ker(s)}$ and define $f^2$ by
$$
f^2\ulam(x,y)=F^2\ulam(x,y)-i(s(F^2\ulam(x,y))).
$$
It is easy to see that $f$ is a $2CL_\infty$-algebra morphism.
Furthermore, $S$ also preserves the identity morphisms and the
composition of morphisms. Thus, $S$ is a functor from {\bf LieC2Alg} to
{\bf CL$_\infty$}.

One show that there are natural isomorphisms $\alpha:T\circ S\Longrightarrow 1_{{\bf LieC2Alg}}$ and $\beta:S\circ
T\Longrightarrow 1_{{\bf 2CL_\infty}}$. This complete the proof.
\qed

\section{Construction of conformal Lie 2-algebras}

In this section, we present concrete examples and special cases of conformal Lie 2-algebras. These include Lie conformal algebras with 3-cocycles, crossed modules of Lie conformal algebras, string Lie conformal algebras, and conformal omni-Lie algebras. Finally, we demonstrate that conformal Lie 2-algebras can be derived from any Leibniz conformal algebra.

\subsection{Skeletal conformal Lie 2-algebras}
In the case of a 2-term conformal $L_{\infty}$-algebra with $\dM=0$, it is referred to as \emph{skeletal}. From conditions $(a)$ and $(g)$, it can be deduced that $\calV_0$ constitutes a Lie conformal algebra. Furthermore, conditions $(b)$ and $(h)$ indicate that $\calV_1$ serves as a representation of $\calV_0$ through the action established by $x \trr \ulam h := l^2_{\lam}(x, h)$. Moreover, condition $(i)$ can be expressed as a 3-cocycle condition in the Lie conformal algebra cohomology of $\calV_0$ with respect to $\calV_1$ as its values.

\begin{Proposition}
Skeletal 2-term conformal $L_{\infty}$-algebras are in one-to-one correspondence with
quadruples $(R, M, \trr\ulam, l^3_{\lam_1,\lam_2})$ where $R$ is a Lie conformal algebra, $M$ is a $\C[\partial]$-module,
$\trr\ulam$ is an action of $R$ on $M$ and $l^3_{\lam_1,\lam_2}$ is a 3-cocycle on $R$ with coefficients in $M$.
\end{Proposition}


Let $M$ be a $\mathbb{C}[\partial]$-module. A conformal bilinear form on $M$ is a $\mathbb{C}$-bilinear map $\langle,\rangle_\lambda: M \times M \rightarrow \mathbb{C}[\lambda]$ such that
$$
\langle\partial v, w\rangle_\lambda=-\lambda\langle v, w\rangle_\lambda=-\langle v, \partial w\rangle_\lambda \quad \text { for all } v, w \in M.
$$
The conformal bilinear form is symmetric if $\langle v, w\rangle_\lambda=\langle w, v\rangle_{-\lambda}$ for all $v, w \in V$.
The conformal bilinear form in a Lie conformal algebra $R$ is called invariant if
$$
\left\langle\left[x_\mu y\right], z\right\rangle_\lambda=\left\langle x,\left[y_{\lambda-\partial} z\right]\right\rangle_\mu
=-\left\langle x,\left[z_{-\lambda} y\right]\right\rangle_\mu
$$
for all $x,y, z\in R$.


\begin{Example}
Given a Lie conformal algebra $(\g,[\cdot _\lambda \cdot])$ with symmetric invariant conformal bilinear form $\langle,\rangle_\lambda$, we construct a 2-term conformal $L_{\infty}$-algebra as follows.
Let $\calV_1=\C, \calV_0=\g, \dM=0$, and define $l^2_{\lam}, l^3_{\lam_1,\lam_2}$ by
\begin{equation}\label{eqn:l2l3string}
l^2_{\lam}(x,y)=[x\ulam y],\quad l^2_{\lam}(x,h)=0,\quad l^3_{\lam_1,\lam_2}(x,y,z)=\langle[x\ulami y],z\rangle\ulamii,
\end{equation}
where $x,y,z\in\g, h\in \C$.
All the conditions in the Definition \ref{2termliealgebra} are satisfied. Thus we obtian a 2-term conformal $L_{\infty}$-algebra $(\C\stackrel{0}{\longrightarrow}\g,l^2_{\lam},l^3_{\lam_1,\lam_2})$. We call it the \emph{string conformal Lie 2-algebra}.
\end{Example}

\subsection{Strict conformal Lie 2-algebras}
Another kind of  2-term conformal $L_{\infty}$-algebra is called \emph{strict} if $l^3_{\lam_1,\lam_2}=0$.
These conformal Lie 2-algebras can be characterized using crossed modules of Lie conformal algebras.

\begin{Definition} Let $(\g,[\cdot,\cdot]_{\g})$ and $(\h,[\cdot,\cdot]_{\h})$ be two Lie conformal algebras.
A crossed module of Lie conformal algebras is a homomorphism of Lie conformal algebras
$\varphi: \h\to \g$ together with an action of $\g$ on $\h$, denoted by $x\trr\ulam h$,
such that
$$\varphi(x\trr\ulam h) = [x\ulam \varphi(h)]_{\g},\quad\varphi(h)\trr\ulam k = [h\ulam k]_{\h},$$
for all $h, k\in\h, x\in\g$.
\end{Definition}

\begin{Proposition} There is an one-to-one correspondence between
strict 2-term conformal $L_{\infty}$-algebras and  crossed modules of Lie conformal algebras.
\end{Proposition}

\pf Let $\calV_1\stackrel{\dM}{\longrightarrow} \calV_0$ be a $2$-term conformal $L_{\infty}$-algebra with $l^3_{\lam_1,\lam_2}=0$. We construct Lie conformal algebras on $\g=\calV_0$ and $\h=\calV_1$ as follows.
The bracket  on $\g$ and $\h$ are defined by
\begin{eqnarray*}
&&[h\ulam k]_{\h}:=l^2_{\lam}(\dM h,k),\quad\forall~x,y\in \h=\calV_1;\\
&&[x\ulam y]_{\g}:=l^2_{\lam}(x,y),\quad\forall~h,k\in \g=\calV_0.
\end{eqnarray*}
By condition $(a)$ and $(g)$ in Definition \ref{2termliealgebra},
it is easy to see that $[\cdot\ulam \cdot]_\g$ satisfies the  Jacobi identity.
By $(h)$, we have
\begin{eqnarray*}
 && [h{}\ulam[k{}\umu l]]-[[h{}\ulam k]{}\ulamu l]- [k{}\umu [h,{}\ulam l]\\
 &=&l^2_{\lam}(\dM h,l^2_{\mu}(\dM k,l))-l^2\ulamu (\dM l^2_{\lam}(\dM h,k),l)-l^2_{\mu}(\dM k,l^2_{\lam}(\dM h,l))\\
 &=&l^2_{\lam}(\dM h,l^2_{\mu}(\dM k,l))-l^2\ulamu (l^2_{\lam}(\dM h,\dM k),l)-l^2_{\mu}(\dM k,l^2_{\lam}(\dM h,l))\\
 &=&0.
\end{eqnarray*}
Thus $[\cdot\ulam\cdot]_\h$ satisfies the  Jacobi identity.
Now let $\varphi=\dM$, then by $(e)$, we have
$$
\varphi([h\ulam k]_\h)=\dM(l^2_{\lam}(\dM h,k))=l^2_{\lam}(\dM h,\dM k)=[\varphi(h)\ulam\varphi(k)]_\g,
$$
which implies that $\varphi$ is a homomorphism of Lie conformal algebras.

Now define the maps of $\trr\ulam:\g\times \h\to \h[\lam]$ by
$$x\trr\ulam h:=l^2_{\lam}(x,h)\in\h,$$
which is an action of $\g$ on $\h$ and it is easy to check that
\begin{eqnarray*}
&&\varphi(x\trr\ulam h)=\dM(l^2_{\lam}(x,h)) = l^2_{\lam}(x, \dM h)= [x\ulam \varphi(h)]_{\g}\\
&&\varphi(h)\trr\ulam k =l^2_{\lam}(\dM h,k)=[h\ulam k]_{\h}.
\end{eqnarray*}
Therefore, we obtain a crossed module of Lie conformal algebras.

Conversely,  a crossed module of Lie conformal algebras
 gives rise to a 2-term conformal $L_{\infty}$-algebra with $\dM=\varphi$,
$\calV_0=\g$ and $\calV_1=\h$, where the brackets are given by
\begin{eqnarray*}
~ l^2_{\lam}(x,y)&:=&[x\ulam y]_{\g},\quad \forall
~x,y\in\g;\\
~l^2_{\lam}(x,h)&:=&x\trr\ulam h,\quad\forall~ x\in\h;\\
~l^2_{\lam}(h,k)&:=&0.
\end{eqnarray*}
The crossed module conditions give various conditions for $2$-term conformal $L_\infty$-algebras with $l^3_{\lam_1,\lam_2}=0$.

\qed

\medskip

Let $\g$ be a Lie conformal algebra, recall that a map $D\in \operatorname{Cend}(\g)$ is called
a conformal derivation if
$$D\ulam[x\umu y]=[(D\ulam x)\ulamu y]+[x\umu D\ulam y].$$
For example, the adjoint action $(\ad x)\ulam$ defined by
$$(\ad x)\ulam(y)=[x\ulam y]$$
from $\g$ to itself is a conformal derivation. This is called the inner conformal derivation.

Denote by $\Der(\g)$ the vector space spanned by all conformal derivations.
We find that $\Der(\g)$ becomes a Lie conformal algebra under the bracket
$$[D\ulam D']=D\ulam D'-D'\uplam D,$$
equivalently,
$$\left[D_\lambda D'\right]_{\lam'}=D_\lambda D'_{{\lam'}-\lambda}-D'_{{\lam'}-\lambda} D_\lambda,$$
where $D,D'$ are conformal derivations of $\g$.

For an operation $\om$ on ${R}$ such that $\om\ulam: {R} \times {R} \to {R}[\lam]$, we define the adjoint operator
$\ad_{\om}: {R}\rightarrow \cendr$ by
$$\ad_{\om\ulam}(x)(y)=\om\ulam(x,y)\in {R}, \quad \forall x,y\in {R}.$$
Then the graph of the adjoint operator
$$\calf_{\om} :=\{\ad_{\om} x+x ~; \, \forall x\in {R}\}\subset \cale:=\cendr\oplus {R}$$
is a subspace of $\cale$.

\begin{Proposition}\label{prop:derivation}
$D$ is a  conformal derivation of ${R}$ if and only if $\calf_{\om}$ is an invariant subspace of $D$, that is $D\circ\calf_{\om}\subseteq\calf_{\om}$, where
the action $D$ on $\calf_{\om}$ is defined by
\begin{equation}\label{daction}
D\ulam\circ (\ad_{\om\umu}(x)+x)=[D\ulam \ad_{\om\umu}(x)]+D\ulam x.
\end{equation}
for all $\ad_{\om\umu}(x)+x\in \calf_{\om\umu}$.
\end{Proposition}

\pf
The right hand side of \eqref{daction}  is belong to $\calf_{\om\umu}$ if and only if
\begin{align*}
[D\ulam \ad_{\om\umu}(x)]&=\ad_{\om\umu}(D\ulam x),
\end{align*}
this is equivalent to
\begin{eqnarray*}
&& D\ulam\ad_{\om\umu}(x)(y)- \ad_{\om\umu}(x)D(y)-\ad_{\om\umu}(D\ulam x)(y)\\
&=&D\ulam[x\umu y]-[x\umu (D\ulam y)]-[(D\ulam x)\ulamu y]\\
&=&0.
\end{eqnarray*}
Thus $D$ is a conformal derivation if and only if $D\circ\calf_{\om\ulam}\subseteq\calf_{\om\ulam}$.
\qed

We call the set of $D\in \operatorname{Cend}(\g)$ such that $D\circ \calf_{\om\ulam}\subseteq \calf_{\om}$ the normalizer of $\calf_{\om}$, which is denoted by $N(\calf_{\om})$.
\begin{Proposition}\label{prop:derivation}
Let $D$ and $D'$ be two conformal derivations. Then $[D\ulam D']$ is also a conformal derivation.
Thus we have $\Der(\g)=N(\calf_{\om})$ is a conformal Lie subalgebra of $\cendr$.
\end{Proposition}
\pf Let $D\circ\calf_{\om}\subseteq\calf_{\om}$ and $D'\circ\calf_{\om}\subseteq\calf_{\om}$, then
\begin{align*}
[D\ulam\ad_{\om\umu}(x)]=\ad_{\om\umu}(D\ulam x),\quad [D'_{\lam'} \ad_{\om\umu}(x)]=\ad_{\om\umu}(D'_{\lam'} x).
\end{align*}
By the  Jacobi identity of $\cendr$, we have
\begin{eqnarray*}
&&[[D\ulam D']_{\lam+\lam'}\ad_{\om\umu}(x)]\\
&=&[D\ulam [D'_{\lam'}\ad_{\om\umu}(x)]]-[D'_{\lam'}[D\ulam\ad_{\om\umu}(x)]\\
&=&[D\ulam \ad_{\om\umu}(D'_{\lam'} x)]-[D'_{\lam'}\ad_{\om\umu}(D\ulam x)]\\
&=&\ad_{\om\umu}(D\ulam D'_{\lam'} x)-\ad_{\om\umu}(D'_{\lam'}D\ulam x)\\
&=&\ad_{\om\umu}([D\ulam D']_{\lam+\lam'} x).
\end{eqnarray*}
This is equivalent to $[D\ulam D']\circ\calf_{\om}\subseteq\calf_{\om}$, so the bracket in $\Der(\g)$ is closed.
Thus $\Der(\g)$ is a Lie conformal algebra as a conformal Lie subalgebra of $\operatorname{Cend}(\g)$.
\qed

\begin{Proposition}
Let $\g$ be a Lie conformal algebra, $\Der(\g)$ and  $\Inn(\g)$ be the set of their conformal derivations and inner conformal derivations.
Then we obtain a crossed module $i:\Inn(\g)\to\Der(\g)$, with $\Der(\g)$ acting  $\Inn(\g)$ by $D\trr\ulam\ad_x=\ad_{D\ulam x}$.
\end{Proposition}

\subsection{Conformal omni-Lie algebras}
The notion of omni-Lie algebras was generalized to omni-Lie superalgebras in \cite{ZL}.
In thus subsection, we introduce the concept of conformal omni-Lie algebras and construct
2-term conformal $L_{\infty}$-algebras from them.


Let $\g$ be a Lie conformal algebra and  $M$ be a left $\g$-module. We define an operation $\circ$ on $\cale:=\g\oplus M$
by
\begin{equation}
  \label{eq:Leibnizbracket}
   ({x+u})\circ\ulam ({y+v}) = [{x\ulam y}]+{x\trr\ulam v},
\end{equation}
for all $x, y\in\g$ and $u, v\in M$.
Then we have
\begin{Proposition}
$(\cale,\circ\ulam)$ is a Leibniz conformal algebra.
\end{Proposition}

\pf We have to check the left  Leibniz identity for the operation $\circ\ulam$.
Let $e_1={x+u}, e_2={y+v}, e_3={z+w}$. Then we have
\begin{eqnarray*}
&&\{e_1\circ\ulam  e_2\} \circ\ulamu e_3-e_1\circ\ulam \{ e_2 \circ\umu e_3\}- e_2\circ\umu \{ e_1 \circ\ulam e_3\} \\
&=&([{x\ulam y}]+{x\trr\ulam v})\circ\ulamu ({z+w})-({x+u})\circ\ulam ([y\umu z]+y\trr\umu w)\\
&&- ({y+v})\circ\umu ([x\trr\ulam w]+z\trr\ulam w)\\
&=&[[{x\ulam y}]\ulamu z]-[x\ulam [y\umu z]]-[y\umu [x\ulam z]]\\
&&+[{x\ulam y}]\trr\ulamu  w-x\trr\ulam  (y\trr\umu w)- y\trr\umu (x\trr\ulam w)\\
&=&0.
\end{eqnarray*}
The right hand side of the above equation is zero because $\g$ is a Lie conformal algebra acting on $M$.
\qed

We call this type of Leibniz conformal algebra structure on $\cale$  the \emph{hemisemidirect product} of $\g$ with $M$ as in \cite{KW},
denote it by $\g\ltimes_H M$.
Note that the operation $\circ$ is not  skew-symmetry.
The  skew-symmetrized bracket of it is:
\begin{equation}  \label{eq:bracket}
  \bleft ({x+u})\ulam ({y+v})\bright := [{x\ulam y}]+\half\left({x\trr\ulam v}- {y\trr\uplam u}\right).
\end{equation}
This is called \emph{demisemidirect product} of $\g$ with $M$, denoted by $\g\ltimes_D M$.

For $e_1={x+u}, e_2={y+v}, e_3={z+w}\in \cale$, denote by
\begin{align*}
\tilde{J}(e_1, e_2, e_3):=
\bleft e_1{}\ulami \bleft e_2{}\ulamii e_3\bright\bright-\bleft e_2{}\ulamii \bleft e_1{}\ulami e_3\bright\bright
-\bleft\bleft e_1{}\ulami e_2\bright{}\ulamipii e_3\bright{}.
\end{align*}
Then $\tilde{J}(e_1, e_2, e_3)\neq 0$.
For example, let $e_1=x, e_2=y, e_3=w$, then
$$\tilde{J}(x, y, z)=-\four[x\ulami y]\trr\ulamipii w.$$
Thus  $\g\ltimes_D M$ can not to be a Lie conformal algebra but a Lie conformal 2-algebra, see Theorem \ref{omniare2term} as follows.

Now, for a Lie conformal algebra $\g$ and a left $\g$-module $M$, let
$$\calV_0=\g\ltimes_D M,\quad  \calV_1=M,\quad \dM=i: M\hookrightarrow \g\ltimes_D M$$
where $i$ is the inclusion map  and define
\begin{align*}
l^2_{\lam}=\bleft\cdot\ulam\cdot\bright,\quad l^3_{\lam_1,\lam_2}=\tilde{J}.
\end{align*}

\begin{Theorem}\label{omniare2term}
Let $R$ be a Lie conformal algebra and $M$ a left $\g$-module.   Then we obtain a nontrivial 2-term conformal $L_{\infty}$-algebra $(M\stackrel{i}{\hookrightarrow}\g\ltimes_D M,\, l^2_{\lam}, \,l^3_{\lam_1,\lam_2})$.
\end{Theorem}

\pf
It can be checked that various conditions in Definition \ref{2termliealgebra} hold.
First we have
\begin{eqnarray*}
&&\bleft ({x+u})\ulam ({y+v})\bright+ \bleft ({y+v})\uplam ({x+u})\bright\\
&=& [{x\ulam y}]+\half\left({x\trr\ulam v}- {y\uplam u}\right)+[{y}\uplam {x}]\\
&&+ \half\left({y\uplam u}- {x\trr\ulam v}\right)\\
&= &[{x\ulam y}]+ [{y}\uplam {x}]+\half\left({x\trr\ulam v}- {y\uplam u}\right)\\
&&+\half\left( {y\uplam u}-{x\trr\ulam v}\right)\\
&=&0.
\end{eqnarray*}
Thus condition $(a)$ holds.

Let $e_1={x}, e_2={y}, e_3={z}$ where ${x}, {y}, {z}\in \g$ and $e_4=t\in M$, then we have
\begin{eqnarray*}
&&\bleft {x}\ulami l^3_{\lam_2,\lam_3}({y}, {z},t)\bright -\bleft {y}\ulamii l^3_{\lam_1,\lam_3}({x}, {z},t)\bright \\
&& +\bleft {z}\ulamiii l^3_{\lam_1,\lam_2}({x}, {y},t)\bright+\bleft l^3_{\lam_1,\lam_2}({x}, {y}, {z})\ulamipiipiii  t\bright- l^3_{\lam_1+\lam_2,\lam_3}(\bleft {x}\ulami {y}\bright, {z},t)\\
&&- l^3_{\lam_2,\lam_1+\lam_3}({y},\bleft {x}\ulami {z}\bright,t)-l^3_{\lam_2,\lam_3}({y}, {z},\bleft {x}\ulami t\bright)\\
&&+ l^3_{\lam_1,\lam_2+\lam_3}({x},\bleft {y}\ulamii {z}\bright, t) +l^3_{\lam_1,\lam_3}({x},{z},\bleft {y}\ulamii t\bright)- l^3_{\lam_1,\lam_2}({x}, {y},\bleft {z}\ulamiii t\bright)\\
&=&-\eight {x}\trr\ulami ([{y}\ulamii {z}]\ulamiipiii\trr t) +\eight {y}\trr\ulamii ([{x}\ulami {z}]\ulamipiii\trr t)\\
&&-\eight {z}\trr\ulamiii([{x}, {y}]\trr\ulamipii t)+0 + \four[[{x}\ulami {y}]\ulamipii{z}]\trr\ulamipiipiii t\\
&&- \four[[{x}\ulami{z}]\ulamipiii{y}]\trr\ulamipiipiii t+ \four [[{y}\ulamii{z}]\ulamiipiii {x}]\trr\ulamipiipiii t \\
&&+\eight[{y}\ulamii {z}]\trr\ulamiipiii({x} \trr\ulami t)-\eight [{x}\ulami {z}]\trr\ulamipiii({y}\trr\ulamii t)+\eight [{x}\ulami {y}]\trr\ulamipii ({z}\trr\ulamiii t)\\
&=&-\textstyle{\frac{3}{8}}([{x}\ulami [{y}\ulamii {z}]] - [[{x}\ulami {y}]\ulamipii {z}] -[{y}\ulamii[{x}\ulami {z}]])\trr\ulamipiipiii t\\
&=&0.
\end{eqnarray*}
Thus condition $(i)$ holds.
The other conditions can be checked similarly.
\qed

Furthermore, we define a  symmetric bilinear
form on $\cale$ with values in $M[\lam]$ by
\begin{equation}  \label{eq:symmetric}
  \langle {x+u}, {y+v}\rangle\ulam:=\half({x\trr\ulam v}+ y\trr\ulam u).
\end{equation}
The triple $(\cale, \bleft\cdot\ulam\cdot\bright, \la\cdot,\cdot\ra\ulam)$ is called a {\bf conformal omni-Lie algebra}.

When $\g=\operatorname{Cend}(M)$,  all possible Lie conformal algebra structures on
$M$ can be characterized by means of  the conformal omni-Lie algebra.

For an  operation $\om_\lam: M  \times M \to M[\lam] $, we define the
adjoint operator
$$\ad_{\om\ulam}: M \rightarrow \operatorname{Cend}(M)  ,\quad \ad_{\om\ulam}(x)(y)=\om\ulam(x,y)\in M[\lam] $$
where $x, y\in M$.
 Then the graph of the adjoint operator:
$$\calf_{\om} =\{\ad_{\om} x+x ~; \, \forall x\in M\}\subset \cale = \operatorname{Cend}(M)\op M$$
is  a subspace of  $\cale$.
 Denote $\calf_{\om}^\perp$ the orthogonal complement of $\calf_{\om}$ in $\cale$ with respect to the
 symmetric bilinear form $\la\cdot,\cdot\ra$ on $\cale $ given
in \eqref {eq:symmetric}.

\begin{Proposition}
\label{prop:realize}
With the above notations,
$(M,\om\ulam)$ is a Lie conformal algebra if and only if its graph
$\calf_{\om}$ is maximal isotropic, i.e. $\calf_{\om}=\calf_{\om}^\perp$, and
 is closed with respect to the bracket
$\bleft\cdot\ulam\cdot\bright$.
\end{Proposition}


\pf  If ${\om\ulam}$ is  skew symmetric, i.e. ${\om\ulam}(x,y)+ {\om\uplam}(y,x)=0$, then
\begin{eqnarray*}
\langle\ad_{{\om\ulam}}(x)+x, \ad_{{\om\uplam}}(y)+y\rangle\ulam&=&\half(\ad_{{\om\ulam}}(x)y+ \ad_{{\om\uplam}}(y)x)\\
&=&\half({\om\ulam}(x,y)+ {\om\uplam}(y,x))
\end{eqnarray*}
 This means that  ${\om\ulam}$ is  skew-symmetric if and only if its graph  is isotropic, i.e.
$\calf_{\om}\subseteq\calf_{\om}^\perp$. Moreover, by dimension
analysis, we have $\calf_{\om}$ is  maximal isotropic.

Next let $[x\ulam y] : ={\om\ulam}(x,y)$, we shall check that the  Jacobi
identity on $M$ is satisfied if and only if $\calf_{\om}$ is closed
under bracket (\ref{eq:bracket}) on $\cale $. In fact,
\begin{eqnarray*}
&&\bleft \ad_{{\om\ulam}}(x)+x, \ad_{{\om\umu}}(y)+y\bright\\
&=&[\ad_{{\om\ulam}}(x),\ad_{{\om\umu}}(y)]+\half(\ad_{{\om\ulam}}(x)y- \ad_{{\om\upmu}}(y)x)\\
&=&[\ad_{{\om\ulam}}(x),\ad_{{\om\umu}}(y)]+\half({\om\ulam}(x,y)- {\om\uplam}(y,x))\\
&=&[\ad_{{\om\ulam}}(x),\ad_{{\om\umu}}(y)]+{\om\ulam}(x,y).
\end{eqnarray*}
Thus this bracket is closed if and only if
$$[\ad_{{\om\ulam}}(x),\ad_{{\om\umu}}(y)]=\ad_{{\om\ulamu}}({\om\ulam}(x,y)).$$
In this case, for $\forall z\in V$, we have
\begin{eqnarray*}
&&[\ad_{{\om\ulam}}(x),\ad_{{\om\umu}}(y)](z)-\ad_{{\om\ulamu}}({\om\ulam}(x,y))(z)\\
&=&\ad_{{\om\ulam}}(x)\ad_{{\om\umu}}(y)(z)- \ad_{{\om\umu}}(y)\ad_{{\om\ulam}}(x)(z)-\ad_{{\om\ulamu}}({\om\ulam}(x,y))(z)\\
&=&\ad_{{\om\ulam}}(x){\om\umu}(y,z)- \ad_{{\om\umu}}(y){\om\ulam}(x,z)-{\om\ulamu}({\om\ulam}(x,y),z)\\
&=&{\om\ulam}(x,{\om\umu}(y,z))- {\om\umu}(y,{\om\ulam}(x,z))-{\om\ulamu}({\om\ulam}(x,y),z)\\
&=&[x\ulam [y\umu z]]- [y\umu [x\ulam z]]-[[x\ulam y]\ulamu z]\\
&=&0.
\end{eqnarray*}
This is exactly the  Jacobi identity on $M$.
\qed
\medskip

We define a Dirac structure of $\cendM\op M$ to be any maximal isotropic
subspace $L\subseteq \cendM\op M$ which is closed under the bracket operation.
Then we obtain that $(M,{\om\ulam})$ is a Lie conformal algebra if and only if $\calf_{\om\ulam}$ is a Dirac structure of the conformal omni-Lie algebra $\g\op M$.

For a general characterization for all Dirac structures of $\cale$, we adapt the theory of characteristic pairs developed in \cite{liuDirac}. One can prove the following result.

\begin{Proposition}\label{prop:Dirac3}
There is a one-to-one correspondence between Dirac structures of the
conformal omni-Lie algebra $(\cendM\oplus M, \bleft\cdot\ulam\cdot\bright, \la\cdot,\cdot\ra\ulam)$ and Lie conformal algebra structures on subspaces of $M$.
\end{Proposition}

\subsection{Skew-symmetrization of Leibniz conformal algebras}
From above subsection, we have seen that a conformal omni-Lie algebra $\cale=\g\ltimes_D M$
is in fact the skew-symmetrization of Leibniz conformal algebra $\g\ltimes_H M$ and $M$ is the left center of $\cale$.
This observation prompts the following question: Can we derive a conformal Lie 2-algebra from any given Leibniz conformal algebra?
In the following subsection, we provide a positive answer to this inquiry.

Let $L$ be a Leibniz conformal algebra. We define the \emph{Leibniz kernel} $\Ker(L)$ of $L$ to be the set of elements spanned by
$$\{x\circ\ulam y+ y\circ\uplam x| \forall x, y\in L\}.$$
Note that if $L$ is a Lie conformal algebra, then $\Ker(L)=0$, otherwise if $L$ is a non-Lie Leibniz conformal algebra, then $\Ker(L)\neq 0$.

The \emph{left center}  $Z^l(L)$ of $L$ is defined by
$$Z^l(L) = \{t\in L|t \circ\ulam x = 0, \forall x\in L)\}.$$
It is easy to see that $\Ker(L)$ and $Z^l(L)$ are ideals of $L$ and the quotient algebras $L/\Ker(L)$ and $L/Z^l(L)$ are Lie conformal algebras.
\begin{Proposition}\label{pro:J0}
The Leibniz kernel $\Ker(L)$ is contained in the left center $Z^l(L)$.
Thus for any non-Lie Leibniz conformal algebra, the set $Z^l(L)$ is not empty.
\end{Proposition}
\pf Let $x\circ\ulam y+ y\circ\uplam x\in \Ker(L)$. Then we get
\begin{eqnarray*}
&&\Big(x\circ\ulam y+ y\circ\uplam x\Big) \circ\ulamu z \\
&=& x \circ\ulam (y \circ\umu z)-  y \circ\umu  (x \circ\ulam  z) \\
&&+ y \circ\umu (x \circ\ulam z) -    x \circ\ulam (y \circ\umu z) \\
&=& 0,
\end{eqnarray*}
for all $z\in L$. Thus $x\circ\ulam y+ y\circ\uplam x\in Z^l(L)$. Therefore $\Ker(L)$ is contained in $Z^l(L)$.
\qed
\medskip

For a Leibniz conformal algebra  $(L,\circ\ulam)$, since the operation $\circ\ulam$ is not  skew-symmetry, we introduce the
following  skew-symmetric bracket on $L$ by
\begin{equation}
\Dorfman{x\ulam y}=\half\left(x\circ\ulam y- y\circ\uplam x\right),\quad\forall x,y\in L,
\end{equation}
and let $\tilde{J}$ be given by
\begin{align}
  \tilde{J}(x,y,z):=&\ \  \Dorfman{\Dorfman{x\ulami y}\ulamipii z}-\Dorfman{x\ulami \Dorfman{y\ulamii z}}-\Dorfman{y\ulamii \Dorfman{x\ulami z}}.
\end{align}


Now it is easy to prove that
\begin{Proposition}
If $t\in  Z^l(L)$, i.e. $t\circ x=0$ for all $x\in L$, then we have
\begin{eqnarray*}
\Dorfman{x\ulam t}=\half x\circ\ulam t,\quad \tilde{J}(x,y,t)=-\four(x\circ\ulami y)\circ\ulamipii t.
\end{eqnarray*}
\end{Proposition}

At last, for any non-Lie Leibniz conformal algebra $L$, we construct nontrivial conformal Lie 2-algebras as follows. Let
$$\calV_0=L,\quad \calV_1=\mathrm Z^l(L),\quad \dM=i: \mathrm Z^l(L)\hookrightarrow L,\quad l^2_{\lam}=\bleft\cdot\ulam\cdot\bright,\quad l^3_{\lam_1,\lam_2}=\tilde{J}.$$

\begin{Theorem}\label{thm:main1}
  With the above notations, from a Leibniz conformal algebra  $(L,\circ\ulam)$, we derive a nontrivial 2-term conformal $L_{\infty}$-algebra
   $(Z^l(L)\stackrel{\dM}{\hookrightarrow}L,\, l^2_{\lam}, \,l^3_{\lam_1,\lam_2})$.
\end{Theorem}

\pf By definition of $\dM,~l^2_{\lam}$ and $l^3_{\lam_1,\lam_2}$, it is easy to see that conditions $(a)$--$(h)$ hold.
For condition $(i)$, we  verify the case of $x,y,z\in L$ and $t\in Z^l(L)$ as follows:
\begin{eqnarray*}
   &&\bleft x\ulami l^3_{\lam_2,\lam_3}(y, z, t)\bright -   \bleft y\ulamii l^3_{\lam_1,\lam_3}(x, z,t)\bright\\
   &&+ \bleft z\ulamiii l^3_{\lam_1,\lam_2}(x, y,t)\bright+\bleft l^3_{\lam_1,\lam_2}(x, y, z)\ulamipiipiii t\bright- l^3_{\lam_1+\lam_2,\lam_3}(\bleft x\ulami y\bright, z,t)\\
   &&- l^3_{\lam_2,\lam_1+\lam_3}(y,\bleft x\ulami z\bright, t)- l^3_{\lam_2,\lam_3}(y, z,\bleft x\ulami t\bright)\\
   &&+l^3_{\lam_1,\lam_2+\lam_3}(x,\bleft y\ulamii z\bright, t) +  l^3_{\lam_1,\lam_3}(x, z,\bleft y\ulamii t\bright)- l^3_{\lam_1,\lam_2}(x, y,\bleft z\ulamiii t\bright)\\
   &=&-\four\Big(\bleft x\ulami ((y\circ z)\circ\ulamiipiii t)\bright-   \bleft y\ulamii ((x\circ\ulami z)\circ\ulamipiii t)\bright  \\
   &&+ \bleft z\ulamiii((x\circ\ulami  y)\circ\ulamipii t)\bright-0- (\bleft x\ulami y\bright\circ\ulamipii z)\circ\ulamIII t\\
   &&+  (\bleft x\ulami z\bright\circ\ulamipiii y)\circ\ulamIII t- (y\circ\ulamii z)\circ\ulamiipiii\bleft x\ulami t\bright\\
   && - (x\circ\ulami\bleft y\ulamii z\bright)\circ\ulamIII t +  (x\circ\ulami z)\circ\ulamipiii\bleft y\ulamii t\bright\\
   &&- (x\circ\ulami y)\circ\ulamipii\bleft z\ulamiii t\bright\Big)\\
   &=&-\textstyle{\frac{3}{8}}\Big(x\circ\ulami (y\circ\ulamii z)-(x\circ\ulami y)\circ\ulamipii z- y\circ\ulamii (x\circ\ulami z)\Big)\circ\ulamipiipiii t\\
   &=&0.
\end{eqnarray*}
The other cases can be checked similarly.
The proof is completed.
\qed

\section*{Acknowledgements}
This research was supported by the National Natural Science Foundation of China (No. 11961049).


%

\vskip7pt

\footnotesize{
\noindent College of Mathematics and Information Science,\\
Henan Normal University, Xinxiang 453007, P. R. China;\\
 E-mail address:\texttt{{  zhangtao@htu.edu.cn}}
}

\end{document}